\documentclass[11pt]{article}
\input{epsf.sty}
\usepackage{hyperref}
\usepackage{amsmath, amssymb}
\newtheorem {thm}{Theorem}[section]
\newtheorem {prop}[thm]{Proposition}

\newtheorem {cor}[thm]{Corollary}

\def\Cox{\hfill \Box}

\def\Z{{\Bbb Z}}
\def\R{{\Bbb R}}

\def\E{{\Bbb E}}

\def\0{{\bf 0}}

\def\sb{{\subset}}

\def\d{\delta}

\def\e{\varepsilon}

\def\phi{\varphi}

\def\D{\Delta}

\def\L{\Lambda}

\def\T{\T}

\begin{document}

\title{Continuous interfaces with disorder:\\
Even strong pinning is too weak in $2$ dimensions } 
\author{
Christof K\"ulske
\footnote{
University of Groningen, 
Department of Mathematics and Computing Sciences, 
Blauwborgje 3,   
9747 AC Groningen, 
The Netherlands
\texttt{kuelske@math.rug.nl}, 
\texttt{ http://www.math.rug.nl/$\sim$kuelske/ }}\, and
Enza Orlandi
\footnote{
Dipartimento di Matematica, 
Universitˆ degli Studi "Roma Tre", 
Largo San Leonardo Murialdo, 1, 
00146 Roma,  ITALY, 
\texttt{orlandi@mat.uniroma3.it },
\texttt{http://www.mat.uniroma3.it/users/orlandi/}}
}

\maketitle

\begin{abstract}
We consider statistical mechanics models of continuous height effective interfaces in the presence 
of a delta-pinning of strength $\e$ at height zero. 
There is a detailed mathematical understanding of the depinning transition in $2$ dimensions without 
disorder.    
Then the variance of the interface height w.r.t. the Gibbs measure stays bounded uniformly in the volume for $\e>0$ 
and diverges like $|\log \e|$ for  $\e\downarrow 0$
How does the presence of a quenched disorder term in the Hamiltonian 
modify this transition?  

We show that an arbitarily weak random field term is enough to beat 
an arbitrarily strong delta-pinning in $2$ dimensions and will cause 
delocalization. The proof is based on a rigorous lower bound for the  
overlap between local magnetizations and random fields in finite volume. 
In $2$ dimensions it  implies growth faster than the volume which is a contradiction to localization.    
We also derive a simple complementary inequality which shows that in higher dimensions 
the fraction of pinned sites converges to one with $\e\uparrow \infty$.

\end{abstract}


\section{Introduction} \label{sect:intro}

\subsection{The setup} 
 
The study of lattice effective interface models, continous and discrete, 
has a long tradition in statistical 
mechanics  \cite{Vel06,BrMeFr86, FoLiNi86, Fun, Sheff, BiKo06,BoKu94, BoKu96}.  


The model we study is given in terms of  variables  $\phi_i\in \R$ which, physically speaking,  
are thought to represent
height variables of a random surface at the sites $i\in \Z^d$. 
Mathematically speaking they are just continuous unbounded (spin) variables. 
The  model is defined in terms of: a  pair potential $V$, 
a quenched random term, and a pinning 
term at interface height zero.   

More precisely, we are interested 
in the behavior of the  quenched finite-volume Gibbs measures  
 in a finite volume $\L\sb \Z^d$ with  
fixed boundary condition at height zero, given by
\begin{equation}\begin{split}\label{localspecification}
&\mu_{\e,\L}[\eta](d\phi_{\L})\cr
&=\frac{
e^{ - \frac{1}{4 d }\sum_{\langle i,j \rangle \in \L}V(\phi_i-\phi_j) 
- \frac{1}{4 d} \sum_{  i\in\L,j\in \L^c, |i-j|=1 }V(\phi_i) 
+  \sum_{i\in \L}\eta_i \phi_i}      \prod_{i\in \L}(d\phi_i +\e \d_0(d\phi_i))    }
{Z_{\e,\L}[\eta]}
\cr
\end{split}
\end{equation}
where the partition function $Z_{\e,\L}[\eta]$ denotes the normalization constant
that turns the last expression into a probability measure.  
The Dirac-measures at the interface height zero are multiplied with the parameter 
$\e$, having the meaning of a coupling strength. 
The disorder configuration $\eta=(\eta_i)_{i\in \R^d}$ denotes an arbitrary fixed configuration of external fields, modelling a "quenched" (or frozen) random environment.    

What do we expect for such a model? 
Recall that the variance of a free massless interface in a finite box diverges like the logarithm 
of the sidelength when there are no random fields. Adding an arbitrarily small pinning $\e$ (without disorder) always localizes 
the interface uniformly in the volume, with the variance of the field behaving 
on the scale $|\log \e|$ when $\e$ tends to zero. 
Indeed, there is a beautiful and complete mathematical understanding of the model without disorder, in the case 
of both Gaussian and uniformly elliptic potentials (see \cite{BoVe01, DeuVel00}) with precise asymptotics 
as the pinning force tends to zero. These results follow from the analysis 
of the distribution of pinned sites and the random walk (arising 
from the random walk representation of the covariance of the $\phi_i$'s) 
with killing at these sites.
In this sense there is already a random system that needs to be analyzed 
even without disorder in the original model. 

What do we expect if we turn on randomness in the model and add the $\eta_i$'s ? 
Let us review first what we know about the same model without a pinning force. 
In $d=2$  we recently proved the deterministic lower bound 
$\mu_{\L_N}[\eta](|\phi_0| \geq t \sqrt {\log L} )  \geq c \exp(- c t^2)$ uniformly for any 
fixed disorder configuration $\eta$, for general potentials $V$ (assuming 
not too slow growth at infinity)  \cite{KuOr06}.  So, it is not possible to stabilize 
an interface by cleverly choosing a random field configuration (one could 
think e.g. that this might be possible with a staggered field).  As this result holds 
at any arbitrary fixed configuration here we don't need any assumptions 
on the distribution of random fields. This result clearly excludes the existence of an infinite-volume Gibbs measure describing a two dimensional 
interface in infinite volume in the presence of random fields. 
In another paper  \cite{EnKu06} the question of existence of {\it gradient Gibbs measures}  (Gibbs 
distributions of the increments of the interface) in infinite volume  was raised. 
Note that while interface 
states may not exist in the infinite volume such gradient states may very well exist, 
as the example of the two-dimensional Gaussian free field shows, by computation.  
(For existence beyond the Gaussian case which is far less trivial, see  \cite{Fun,FunSp}.)  
It was proved in \cite{EnKu06} that there are no such gradient Gibbs measures in the random model in dimension $d=2$.

Now, turn to the full model in $d=2$.   
In view of the localization taking place at any positive pinning force $\e$ without disorder, a natural 
guess might be that with disorder at least at very large $\e$ there would be pinning. 
However, we show as  a result of the present paper that  
this is not the case, somewhat to our own surprise, 
and an arbitrarily strong pinning does not suffice to keep the interface bounded.

\subsection{Main results}

\subsubsection*{Delocalization in $d=2$ - superextensivity of the overlap }

Denote by $\L_L$ the square of sidelength $2 L +1$ centered 
at the origin. 

In this subsection we consider the 
disorder average of the overlap in $\L_L$  showing 
that it grows faster than the volume. This in particular implies that in two dimensions there 
is never pinning, for arbitrarily weak random field and arbitrarily large 
pinning forces $\e$. Here is the result. 

\begin{thm}\label{overlapthm2d}  Assume that $\sup_{t}V''(t)\leq 1$, $\liminf_{|t|\uparrow\infty} \frac{\log V(t)}{\log |t|}>1$,   and let 
  $\eta_i$ be symmetrically 
   distributed, i.i.d. 
 with finite second moment. 

Let $d=2$.  Then there is a constant $a>0$, independent of the distribution of the random fields  
and the pinning strength $\e\geq 0$, such that 
\begin{equation}\begin{split}\label{main1}
& \liminf_{L\uparrow \infty } \frac{1}{ L^2 \log L } \sum_{i\in \L_L}
\E \Bigl( \eta_i\mu_{\e,\L_L}[\eta](\phi_i)\Bigr)\geq a \, \E(\eta_0^2).  \cr 
\end{split}
\end{equation}
\end{thm}
\bigskip

Note that the growth condition on $V$ includes the quadratic case and ensures  
the finiteness of the integrals appearing in (\ref{localspecification})
for all arbitrarily fixed choices of $\eta$, even at $\e=0$.  

Generalizations to interactions that are non-nearest neighbor 
are obvious; all results go through e.g. for 
finite range and we skip them in this presentation 
for the sake of simplicity. 
We like to exhibit the case of Gaussian random fields (and not necessarily Gaussian potential $V$) 
since the bound acquires a form that looks even more striking because it becomes 
independent of the size of the variance of the $\eta_i$'s (as long as this is strictly positive). 

\begin{cor}
Let us assume that the random fields $\eta_i$ have an i.i.d.  Gaussian distribution with mean zero and 
strictly positive variance of arbitrary size. 

Then, with the same constant  $a$ as above, we have the bound   
\begin{equation}\begin{split}\label{cor1}
& \liminf_{L\uparrow \infty } \frac{1}{ L^2 \log L } \sum_{i\in \L_L}
 \E \Bigl( \mu_{\e,\L_L}[\eta](\phi^2_i)
-  \mu_{\e,\L_L}[\eta](\phi_i)^2 \Bigr) 
\geq a >0 \cr 
\end{split}
\end{equation}
for any $0\leq \e<\infty$. 
\end{cor}

(\ref{cor1}) follows from (\ref{main1}) by partial integration w.r.t. the Gaussian disorder 
average (transforming the overlap into the variance of the $\phi_i$'s).  

Note that, even in the unpinned case of $\e=0$,  Theorem \ref{overlapthm2d} 
is not entirely trivial in the case of general potentials $V$. Here it provides 
an alternative simple way to see the delocalization in the presence of 
random fields (while the explicit lower bound on the tails of \cite{KuOr06} provides
more information.)

\bigskip 

\subsubsection*{Lower bound on overlap in $d\geq 3$ }

The analogue of Theorem \ref{overlapthm2d}  for higher dimensions
is the following.

\begin{thm}\label{overlapthmdgeq3} Let $d\geq 3$ and let $\e\geq 0$ be arbitrary 
and assume the same conditions on $V$ and $\eta_i$ as in Theorem \ref{overlapthm2d}.  

Then there are positive constants $B_1, B_2 <\infty $,  
independent of the distribution of the random fields  
and the pinning strength $\e\geq 0$, such that 
\begin{equation}\begin{split}\label{main3}
& \liminf_{L\uparrow \infty } \frac{1}{ L^d } \sum_{i\in \L_L}
\E \Bigl( \eta_i\mu_{\e,\L_L}[\eta](\phi_i)\Bigr)\geq  \frac{\E(\eta_0^2) \, (-\D^{-1})_{0,0}}{2}
- \log (B_1 + B_2 \e)  \cr 
\end{split}
\end{equation}
where the positive constant 
$(-\D^{-1})_{0,0}$ is the diagonal element of the inverse of the infinite-volume lattice Laplace 
operator whose existence is guaranteed in  $d\geq 3$. 

\end{thm}

\subsubsection*{Lower bound on the pinned volume in $d\geq 3$ }

We complement the previous lower bounds on the overlaps 
which are depinning-type of results by a pinning-type result.  
It is a lower bound on the disorder average 
of the quenched Gibbs-expectation of the 
fraction of pinned sites. While we needed an upper bound 
on the interaction potential $V$ before we are assuming now a lower bound on 
$V$.

\begin{thm}\label{345}  Let $d\geq 3$. Assume that $\inf_{t}V''(t)= c_-  >0$ 
 and let   $\eta_i$ be 
 symmetrically  distributed, i.i.d. 
 with finite second moment. 

Then there exist dimension-dependent constants $C_1,C_2>0$, independent 
of the distribution of the disorder,  such that, for all $\e$ and for all volumes $\L$, 
the disorder average of the fraction of pinned sites obeys the estimate 
\begin{equation}\begin{split}\label{allyouneed}
&\frac{1}{|\L|}\sum_{i\in \L}\E\Bigl( \, \mu_{\e,\L}[\eta](\phi_i=0) \Bigr) 
\geq 1- \frac{C_1 + C_2 \E (\eta_0^2) }{\log \e}. 
\end{split}
\end{equation}

\end{thm}

This shows pinning for the large $\e$ regime in the "thermodynamic sense" that the 
fraction of pinned sites can be made arbitrarily close to one, uniformly 
in the volume. As usual this result does not allow to make statement 
about the Gibbs measure itself.

\bigskip

The proofs follows from "thermodynamic reasoning". 
The first "depinning-type" result follows from taking the log of the partition function 
and differentiating and integrating back w.r.t. the coupling strength of the random fields. 
Exploiting the linear form of the random fields, convexity, comparison 
of non-Gaussian with the Gaussian partition functions, and asymptotics of Green's functions
the results follow, see Chapter 2.

\bigskip
\bigskip
\bigskip

\section{Proof of Depinning-type results}

The estimates in formulas (\ref{main1}), (\ref{cor1}), and (\ref{main3})
are immediate consequences of the following fixed-disorder estimate.

\begin{prop} For any dimension 
$d$,  there are  constants $C_{\text{nG}, d}<\infty $ and $c_{\text{G},d}>0$
such that, for 
all fixed configurations of local fields $\eta$, we have  
\begin{equation}\begin{split}\label{localspecification1}
& \frac{1}{2}\sum_{i,j\in \L}(-\D_{\L})^{-1}_{i,j}\eta_i\eta_j - |\L| \log\frac{ C_{\text{nG}, d} + \e}{c_{\text{G},d}}
\leq \sum_{i\in \L} \eta_i\mu_{\e,\L}[\eta](\phi_i). 
\end{split}
\end{equation}

\end{prop}

\bigskip 

{\bf Proof of the Proposition: }
Let us see  what comes out when we differentiate and integrate 
back the free energy in finite volume w.r.t. strength of the random fields. 

 \begin{equation}\begin{split}\label{q1}
& \frac{d}{dh}\log Z_{\e,\L}[h\eta]
=\sum_{i\in \L}\eta_i\mu_{\e,\L}[h\eta](\phi_i).
\end{split}
\end{equation}
At every fixed $\eta$, this quantity is a monotone function 
of $h$, which is seen by another differentiation w.r.t. $h$ 
which produces the variance. 
We have 
 \begin{equation}\begin{split}\label{q3}
& \log \frac{Z_{\e,\L}[\eta]}{Z_{\e,\L}[0]}
=\sum_{i\in \L}\int_0^1 dh\eta_i\mu_{\e,\L}[h \eta](\phi_i)
\leq \sum_{i\in \L}  \eta_i \, \mu_{\e,\L}[\eta](\phi_i).\cr
\end{split}
\end{equation}
We note the lower bound on the numerator which we get by dropping the pinning term, giving us 
 \begin{equation}\begin{split}\label{q4}
& Z_{\e,\L}[\eta] 
\geq Z_{\e=0,\L}[\eta] \cr
&\geq Z_{\e=0,\L}^{\text{Gauss}}[\eta] \cr
&=\exp\Bigl( \frac{1}{2}\sum_{i,j\in \L}(-\D_{\L})^{-1}_{i,j}\eta_i\eta_j 
\Bigr)Z_{\e=0,\L}^{\text{Gauss}}[0] \cr
&\geq \exp\Bigl( \frac{1}{2}\sum_{i,j\in \L}(-\D_{\L})^{-1}_{i,j}\eta_i\eta_j 
\Bigr)c_{\text{G},d}^{|\L|}\cr
\end{split}
\end{equation}
Here we have denoted by $ Z_{\e=0,\L}^{\text{Gauss}}[\eta] $ 
the Gaussian partition function with potential $V(t)=\frac{t^2}{2}$. 

Further  we used that the lower bound 
on $V(t)$ taken from the hypothesis implies that, for any partition function 
in any volume $D$, we have $ Z_{\e=0,D}[0] \leq C_{\text{nG},d}^{|D|}$. 
This gives 
\begin{equation}\begin{split}\label{q5}
&Z_{\e,\L}[0] 
= \sum_{A\subset \L} \e^{|A|} Z_{\e,\L\backslash A}[0] \cr
&\leq   \sum_{A\subset \L} \e^{|A|} C_{\text{nG},d}^{|\L\backslash A|}
=(C_{\text{nG},d}+\e )^{|\L|}. 
\end{split}
\end{equation}
So the desired estimate on the overlap follows from (\ref{q3}),(\ref{q4}),(\ref{q5}).
This concludes the proof of the Proposition. 
$\Cox$

\bigskip 
\bigskip 
It is easy to obtain the Theorems
\ref{overlapthm2d} and \ref{overlapthmdgeq3}  
from the proposition. Indeed, 
taking a disorder average  we have 
\begin{equation}\begin{split}\label{localspecification2}
& \frac{\E (\eta_0^2)} {2}\sum_{i\in \L}(-\D_{\L})^{-1}_{i,i} -  |\L|\log\frac{ C_{\text{nG}, d} + \e}{c_{\text{G},d}}
\leq \sum_{i\in \L} \E \Bigl( \eta_i\mu_{\e,\L}[\eta](\phi_i)\Bigr). 
\end{split}
\end{equation}
Now use the asymptotics of the Green's-function 
in a square $(-\D_{\L_L})^{-1}_{i,i}\sim \log L$ at fixed $i$ 
to get the first theorem. 
The proof of the case $d\geq 3$ follows from the existence of the infinite-volume Green's-function 
in $d\geq 3$.  

\bigskip 

Finally let us note in passing that a constant magnetic field 
is always winning against an arbitrarily strong pinning, and even more strongly 
than a random field.  
Indeeed, let $d\geq 2$, let $\eta_i=h\geq 0$ for all sites $i$ 
and let $\e\geq 0$ be arbitrary.   
Then, there is a constant $c_d >0$, independent of $h$ and $\e$, such that 
\begin{equation}\begin{split}\label{localspecification3}
& \liminf_{L\uparrow \infty } \frac{1}{ L^{d+2} } \sum_{i\in \L_L}
\mu_{\e,\L_L}[h](\phi_i)\geq  c_d h.
 \cr 
\end{split}
\end{equation}
This again follows from the Proposition, using $\sum_{i,j\in \L}(-\D_{\L_L})^{-1}_{i,j} \sim L^{d+2}$.

\bigskip 
\bigskip

\section{Proof of Pinning-type results } 

To prove the lower bound on the fraction of pinning sites 
in dimension $d\geq 3$  given in Theorem \ref{345} we will in fact prove the following 
fixed-disorder lower bound:

For all finite volumes $\L$ and 
for all realizations $\eta$ we have, for any $\e_0>0$  
 \begin{equation}\begin{split}\label{laballballon}
&\frac{1}{|\L|}
\sum_{i\in \L}\mu_{\e,\L}[\eta](\phi_i=0)\cr
&\geq \frac{1}{\log \frac{\e}{\e_0}}
\Biggl(\log \frac{\e  c_-^\frac{d}{2} }{ \Bigl(1+  \frac{\e_0 c_-^\frac{d}{2}}{(2\pi)^\frac{d}{2}} \Bigr) C_{\text{G},d}}
- \frac{1}{2 c_- |\L|}\sum_{i,j\in \L}(-\D_{\L})^{-1}_{i,j}\eta_i\eta_j  \Biggr) 
\end{split}
\end{equation}
with a constant $C_{\text{G},d}$ defined in (\ref{18}). 

Taking a disorder-expectation  
(\ref{allyouneed}) follows by the finiteness of Green's function in 
the infinite volume $(-\D_{\Z^d})^{-1}_{0,0}$ with $\e_0=1$. 
$\Cox$

\bigskip 

{\bf Proof of (\ref{laballballon}): } 
The proof is based on the trick to differentiate and integrate back the log of the partition 
function, now w.r.t. $\e$: Differentiation gives 
 \begin{equation}\begin{split}\label{localspecification4}
&\e \frac{d}{d\e}\log Z_{\e,\L}[\eta]
= \sum_{i\in \L}\mu_{\e,\L}[\eta](\phi_i=0).
\end{split}
\end{equation}
We integrate this relation back, and it will be important 
for us to do it starting from a positive $\e_0>0$. So we get     
 \begin{equation}\begin{split}\label{12}
&\log \frac{Z_{\e,\L}[\eta]}{ 
Z_{\e_0,\L}[\eta]}
= \int_{\e_0}^\e \frac{d \tilde \e}{\tilde \e}\sum_{i\in \L}\mu_{\tilde \e,\L}[\eta](\phi_i=0)
\leq \log \frac{\e}{\e_0}\cdot \sum_{i\in \L}\mu_{\e,\L}[\eta](\phi_i=0)
\end{split}
\end{equation}
where we have used that $\sum_{i\in \L}\mu_{\tilde \e,\L}[\eta](\phi_i=0)$ is a monotone 
function of $\tilde \e$. 
Note that  the integrand itself is not a monotone function. 
(Compare \cite{CaVel00} for a related non-random pinning scenario, with back-integration from zero.)

Now we have the trivial lower bound obtained by keeping 
only the contribution in the expansion where all  sites are pinned, i.e. 
  \begin{equation}\begin{split}\label{13}
&Z_{\e,\L}[\eta]\geq \e^{|\L|}.\cr
\end{split}
\end{equation}
For the upper bound on the partition function 
of the full model (at $\e_0$) we first
use the lower bound on the potential $V(t)\geq \frac{c_{-} t^2}{2}$ 
giving us a comparison with a Gaussian partition function with curvature $c_-$: 
  \begin{equation}\begin{split}\label{14}
&Z_{\e_0,\L}[\eta]
\leq Z_{\e_0,\L}^{\text{Gauss},c_-  }[\eta].\cr
\end{split}
\end{equation}
It is a simple matter to rescale the Gaussian curvature away
 \begin{equation}\begin{split}\label{15}
& Z_{\e_0,\L}^{\text{Gauss},c_-  }[\eta]
={c_-}^{-\frac{d |\L|}{2}} Z_{\e_0 {c_-}^{\frac{d}{2}} ,\L}^{\text{Gauss}}[{c_-}^{-\frac{1}{2}} \eta]
\cr
\end{split}
\end{equation}
where the partition function on the r.h.s. is taken with unity curvature potential. 
For the Gaussian partition function we claim the upper bound (writing 
again in the original parameters) of the form 
\begin{equation}\begin{split}\label{16}
&Z_{\e ,\L}^{\text{Gauss}}[\eta] \leq \Bigl(1+  \frac{\e}{(2\pi)^\frac{d}{2}} \Bigr) ^{|\L|} 
Z_{\e=0,\L}^{\text{Gauss}}[\eta].  \cr
\end{split}
\end{equation}
Here is an elementary proof: 
We will replace successively the single-site integrations 
involving the Dirac measure by integrations only over the Lebesgue measure with 
the appropriately adjusted prefactor. Indeed, 
consider one site $i$ and compute the contribution 
to the partition function while 
fixing the values of $\phi_j$ for $j$ not equal to $i$. 
Then use that  
\begin{equation}\begin{split}\label{17}
&\int\Bigl(  d\phi_i + \e \d_0(d\phi_i)
\Bigr)\exp\Bigl( -\frac{\phi_i^2}{2}+ (\sum_{j\sim i}\phi_j + \eta_i )\phi_i 
\Bigr) \cr
&=(2\pi)^\frac{d}{2}\exp\Bigl( \frac{(\sum_{j\sim i}\phi_j + \eta_i )^2}{2}
\Bigr) + \e 
 \cr
 &\leq \Bigl(1+  \frac{\e}{(2\pi)^\frac{d}{2}} \Bigr) 
 (2\pi)^\frac{d}{2}\exp\Bigl( \frac{(\sum_{j\sim i}\phi_j + \eta_i )^2}{2}\Bigr)\cr
 &=  \Bigl(1+  \frac{\e}{(2\pi)^\frac{d}{2}} \Bigr) 
  \int d\phi_i  
\exp\Bigl( -\frac{\phi_i^2}{2}+ (\sum_{j\sim i}\phi_j + \eta_i )\phi_i 
\Bigr) \cr
\end{split}
\end{equation}
and iterate over the sites.

For the Gaussian unpinned partition function use 

\begin{equation}\begin{split}\label{18}
&Z_{\e=0,\L}^{\text{Gauss}}[\eta]  
=\exp\Bigl( \frac{1}{2}\sum_{i,j\in \L}(-\D_{\L})^{-1}_{i,j}\eta_i\eta_j 
\Bigr)Z_{\e=0,\L}^{\text{Gauss}}[0] \cr
&\leq \exp\Bigl( \frac{1}{2}\sum_{i,j\in \L}(-\D_{\L})^{-1}_{i,j}\eta_i\eta_j 
\Bigr)C_{\text{G},d}^{|\L|}\cr
\end{split}
\end{equation} 
with a suitable constant. 
From here (\ref{allyouneed}) follows from (\ref{12},\ref{13},\ref{14},\ref{15},\ref{16},\ref{18})
$\Cox$ 

\bigskip 
\bigskip

{\bf Acknowledgements: } 
The authors thank Pietro Caputo for an interesting discussion and Aernout van Enter 
for comments on a previous draft of the manuscript.  C.K. thanks 
the university Roma Tre for hospitality.

\end{document}